\documentclass[a4paper,12pt]{article}

\newtheorem{theorem}{Theorem}
\newtheorem{lemma}{Lemma}

\newtheorem{corollary}{Corollary}

\newtheorem{remark}{Remark}

\def\C{\textup{\hspace*{.2em}\mbox{l\hspace*{-.5em}C$\,$}}}
\def\R{\textup{\hspace*{.2em}\mbox{I\hspace*{-.1em}R$\,$}}}
\newcommand{\cD}{{\mathcal D}}

\newcommand{\cX}{{\mathcal X}}

\usepackage{amsmath}
\usepackage{amssymb}
\usepackage{amsfonts}
\usepackage{graphicx}
\usepackage{mathrsfs}
\usepackage{color}
\usepackage{ifpdf}

\author{Olivier Bachelier$^1$, Didier Henrion$^{2,3,4}$, Nima Yeganefar$^1$, Driss Mehdi$^1$}

\title{On the solutions to complex parameter-dependent LMIs involved in the stability analysis\\ of 2D discrete models} 

\begin{document}

\date{Version of \today}
\maketitle

\footnotetext[1]{University of Poitiers, LIAS-ENSIP, B\^atiment B25, 2 rue Pierre Brousse, TSA 41105, 86073 Poitiers cedex, France.}  
\footnotetext[2]{CNRS, LAAS, 7 avenue du colonel Roche, F-31400 Toulouse, France.}
\footnotetext[3]{Universit\'e de Toulouse; F-31400 Toulouse; France.}
\footnotetext[4]{Faculty of Electrical Engineering, Czech Technical University in Prague,
Technick\'a 2, CZ-166 26 Prague, Czech Republic.}

\begin{abstract}               
The aim of this short communiqu\'e is to adapt a result established by Bliman,
related to the possible approximation of the solutions to real-parameter-dependent linear matrix inequalities (LMIs), to the
special context of stability analysis of 2D discrete Roesser models.
While Bliman considered the case of LMIs involving several {\em real} parameters, which is especially crucial for the analysis
of linear systems against parametric deflections, the stability of Roesser models leads to consider
LMIs with only one single {\em complex} parameter.
Extending the results from real parameters to complex ones is not straightforward in our opinion. 
This is why the present note discusses precautions to be taken concerning this case
before applying the results in a 2D context. Actually, it is shown that a well-known condition for
{\em structural} stability of a 2D discrete Roesser can be relaxed into an LMI system
whose solution polynomially depends on a single complex parameter over the unit circle.
\end{abstract}

\section{Motivation and introduction}
\label{sec:motiv}

Multidimensional models raised a growing interest within the community of automatic control
during the last three decades, due to the wide field of their applications (signal and image processing, seismic
phenomena,~... see~\cite{Kac85,gw01}). All systems for which information propagates in more than one direction
are eligible to be described by such models. Among them, the 2D models are probably the most
studied ones, since they can be used to describe the behaviours of, {\em e.g.} Iterative Learning Control (ILC)-schemes
or so-called repetitive processes~\cite{gw01,Rogers2007-1}. In the present contribution,
we focus on 2D Roesser models. These models comply with

\begin{equation}
\label{eq:Roesser}
\left[\begin{array}{c}
       {\bf q}_1x_1(j_1,j_2)\\
{\bf q}_2x_2(j_1,j_2)
      \end{array}\right]=\left[\begin{array}{cc}
A_{11} & A_{12}\\
A_{21} & A_{22}
\end{array}\right]\left[\begin{array}{c}
       x_1(j_1,j_2)\\
x_2(j_1,j_2)
      \end{array}\right],
\end{equation}
where $x_i(j_1,j_2)\in{\R}^{k_i}$ is the state subvector associated with dimension $i$.
${\bf q}_i$ is either a derivative operator
or a shift operator along the $i$th dimension.
Therefore,~(\ref{eq:Roesser}) can be continuous, discrete or mixed.
The state vector is divided into two subvectors corresponding to the propagation
of the information in two directions. The matrices $A_{ij}$ are here assumed
to be real.\\
Besides, we define the following subsets of the extended complex plane $\overline{\C}=\C\cup\{\infty\}$,
\begin{equation}
\label{eq:defDi}
\left\{\begin{array}{lll}
 \cD_i & = & \{\lambda\in\overline{\C}: f_i(R_i,\lambda)<0\},\\
 \cD_i^C & = & \{\lambda\in\overline{\C}: f_i(R_i,\lambda)\geq0\},\\
 \partial\cD_i & = & \{\lambda\in\overline{\C}: f_i(R_i,\lambda)=f_i(\hat{R}_i,\lambda)=0\},\\
 \cD_i^{\diamond} & = & \{\delta\in\overline{\C}: \delta^{-1}\in\cD_i^{C}\}= \{\delta\in\overline{\C}: f_i(\hat{R}_i,\delta)\geq 0\}
 \end{array}\right.
\end{equation}
with
\begin{equation}
 f_i(R_i,\lambda)=\left[\begin{array}{c}
                       \lambda\\
                       1
                      \end{array}\right]'R_i\left[\begin{array}{c}
                       \lambda\\
                       1
                      \end{array}\right],\, R_i=\left[\begin{array}{cc}
                      r_{i_{11}} & r_{i_{10}}\\
                      r_{i_{10}} & r_{i_{00}}
                      \end{array}\right]\in\R^{2\times 2},\,\hat{R}_i=\left[\begin{array}{cc}
                      r_{i_{00}} & r_{i_{10}}\\
                      r_{i_{10}} & r_{i_{11}}
                      \end{array}\right].
\end{equation}
The matrices $R_i$ are restricted to cases where either $R_i=\left[\begin{array}{cc} 0 & 1\\ 1 & 0\end{array}\right]$
or $R_i=\left[\begin{array}{cc} 1 & 0\\0 & -1\end{array}\right]$. In such cases $\cD_i$ corresponds to the usual stability region,
{\em i.e.} to the open left half complex plane or the open unit disc respectively. $\cD_i^C$ is the complementary region,
$\partial\cD_i\subset\cD_i^C$
is the boundary of $\cD_i$ and $\cD_i^{\diamond}$ is a set described by $\lambda^{-1}$ when $\lambda$ describes $\cD_i^C$.\\
In order to make sense, to each operator ${\bf q}_i$ we associate a region $\cD_i$ as follows:

\begin{equation}
\label{eq:choices_q}
  \left\{\begin{array}{lllllll}
   {\bf q}_ix_i(j_1,j_2) & = & \displaystyle\frac{\partial}{\partial j_i}x_i(j_1,j_2) & \Rightarrow &
   \cD_i & = & \overline{\C}^-,\,\, i=1,2,\\
   {\bf q}_1x_1(j_1,j_2) & = & x_1(j_1+1,j_2) & \Rightarrow &																	
   \cD_1 & = & {\bf D},\\
   {\bf q}_2x_2(j_1,j_2) & = & x_1(j_1,j_2+1) & \Rightarrow &
   \cD_2 & = & {\bf D},
         \end{array}\right.
\end{equation}
where $\overline{\C}^-$ denotes the achieved open left half complex plane and ${\bf D}$ is the unit disc.
Hence, the continuous, discrete and mixed cases can all be taken into account in this framework.
A frequently used stability criterion is as follows.

\begin{lemma}\label{thm:2Dclassic}
System~(\ref{eq:Roesser}) is stable if and only if
$A_{22}$ is $\cD_2$-stable and
$M(\delta)=A_{11}+A_{12}(I - \delta A_{22})^{-1}\delta A_{21}$ is $\cD_1$-stable
for any $\delta\in\partial\cD_i$ (or any $\delta\in\cD_i^{\diamond}$.
\end{lemma}
The justification can be found for example in~\cite{Huang72} for the countinous case, and in~\cite{siljak75}
for the other cases. From this let us derive the next lemma.

\begin{lemma}
\label{th:intro}
System~(\ref{eq:Roesser}) is stable if and only if there exist a positive definite matrix
$Y$ and a positive definite parameter-dependent Hermitian matrix $P(\delta)$ such that

\begin{equation}
\label{eq:steinA22}
r_{2_{00}}Y+r_{2_{10}}A'_{22}Y+r_{2_{10}}YA_{22}+r_{2_{11}}A'_{22}YA_{22}<0
\end{equation}
and $\forall\delta\in\partial\cD_2 (resp. \forall \delta\in\cD_2^{\diamond})$,
{\small
\begin{equation}
\label{eq:stein}
r_{1_{00}}P(\delta)+r_{1_{10}}M'(\delta)P(\delta)+r_{1_{10}}P(\delta)M(\delta)+r_{1_{11}}M'(\delta)P(\delta)M(\delta)<0.
\end{equation}}
\end{lemma}

{\bf Proof.} This is just an interpretation of Lemma~\ref{thm:2Dclassic} in terms of Lyapunov's or Stein's inequalities.~$\Box$

~\\
Lemma~\ref{th:intro} is interesting but involves condition~(\ref{eq:stein}) which is not computationally
tractable due to its dependence on the complex parameter $\delta$. 
However, if one could prove that $P(\delta)$ might comply to a specified non conservative structure,
it could allow the research for tractable relaxations.
The establishment of such relaxations is the underlying motivation behind this work. But the present note focuses
on the preliminary step {\em i.e.} the relevance of the structure that can be assumed for $P(\delta)$. 
In the sequel we discuss this issue in order to prove the main result which is now stated.
\begin{theorem}
\label{th:th2D}
System~(\ref{eq:Roesser}) is stable if there exist positive definite matrices
$Y$ and

\begin{equation}
\label{eq:Pdedelta}
P(\delta)=\displaystyle\left(\displaystyle\sum_{i=1}^{\nu}P_i\delta^i\right)^H,\, P_i\in\C^{n\times n},
\end{equation}
(where $X^H$ denotes $X+X'$) such that~(\ref{eq:steinA22}) and~(\ref{eq:stein})
are satisfied.
Moreover, the condition is also necessary for a sufficiently large value of the degree $\nu$.
\end{theorem}

Theorem~\ref{th:th2D} was used in the discrete case in \cite{nDS13} to derive a simple and rather tractable hierarchy of sufficient LMI conditions
for (\ref{eq:steinA22}-\ref{eq:Pdedelta}), which tends to necessity, by using a classical S-procedure as formulated in~\cite{scherer01}.
Unfortunately, in \cite{nDS13}, Theorem~\ref{th:th2D} was taken at a face value with no actual proof.
More precisely, the result in~\cite{nDS13} was based upon~\cite{bliman04} but we realized that it was not
technically correct. This paper aims at filling this gap.
To fully motivate Theorem~\ref{th:th2D}, we have to explain the interest in applying the S-procedure whereas there exist more
sophisticated relaxation schemes such as moment-SOS (sum of squares), or Lasserre's hierarchy~\cite{Lasserre2009}. Indeed, it could be possible
to write $\delta=a+{\bf i}b$ with $(a;b)\in\R^2$ and {\bf i} being the imaginary unit, and then to prove the existence of a solution
to (\ref{eq:stein}) that would be polynomial with respecto to $a$ and $b$ (see for instance \cite{Scheiderer2006}).
Such a polynomial could be computed by using SOS techniques.
But the interest in using the S-procedure associated with the Linear Fractional Representation (LFR) framework to eliminate $\delta\in\C$
from the condition is that this relaxation not only possibly leads to non conservative LMIs in the problem that we aim at tackling, but, moreover,
straightforwardly yields
clear expressions of these LMIs~\cite{nDS13}.
The derivation of these expressions is a key step in the design of stabilizing control laws as highlighted in~\cite{nDS15}
where a necessary and sufficient LMI condition for state feedback stabilization is obtained... provided that Theorem~\ref{th:th2D}
is proved! Hence our motivation to clarify this technical aspect.
\\
The LMIs are a tool whose utility and popularity are no longer to be demonstrated~\cite{bgfb94}.
They consitute an effective approach to solve numerous problems in automatic control. Some of those problems unfortunaltely
involve parameter-dependent LMIs. Typically, analyzing the stability of linear systems against parameter uncertainties amounts to assessing the existence
of a solution ({\em i.e.} a Lyapunov function) which itself depends on these parameters~\cite{fag96,goh98,daaber01},
which is a difficult problem. In order to efficiently relax such a problem, the question of the way this solution
could depend on the parameters is of course important. Bliman made a very important
step in the understanding of such a question. Indeed, in~\cite{bliman04}, he proved that when
the LMI depends on real parameters, its solution can be considered polynomial with respect
to the parameters without conservatism. The question of its minimal degree is still open.
In the remainder of the paper, inspired by the work of Bliman~\cite{bliman04} and motivated by our will to fully justify
the work in~\cite{nDS13}, we discuss the possibility
for the solutions to complex parameter-dependent LMIs
to be assumed polynomial with respect to the complex parameters.
The purpose is not to provide a complete complex counterpart with a complete proof, which would be very redundant with Bliman's work,
but to insist on the steps where some adaptations are required and to highlight the currently existing limits
which prevent the complete generalization of Bliman's result.
But it must be kept in mind that
the case with only one complex parameter
is of special interest for us to prove Theorem~\ref{th:th2D}.

\section{Bliman's result}

In this section, we recall Bliman's theorem and give a very brief outline of his proof.
Let us define the following expression:

{\scriptsize
\begin{equation}
\label{eq:LMIreelle}
G(x,\delta)=G_0(\delta)+\displaystyle\sum_{i=1}^p x_iG_i(\delta),\delta\in{K},x=\left[\begin{array}{c}
                                 x_1\\
\vdots\\
x_p
\end{array}\right]\in\R^p,
\end{equation}}
where $K$ is a compact subset of $\R^m$ and $G_0$, $G_1$, ..., $G_p$ are mappings
defined in $K$ and taking values in the set of symmetric matrices of $\R^{n\times n}$.
\begin{theorem}
\label{th:Bliman}
\cite{bliman04}
Assume that mappings $G_i$, $i=1$, ..., $p$, are continuous. If for any $\delta\in K$,
there exists $x(\delta)$
such that $G(x(\delta),\delta)>0$, then there exists a polynomial function
$x^*:K\rightarrow \R^p$ such that
for any $\delta\in K$, $G(x^*(\delta),\delta)>0$.
\end{theorem}
The result is very strong for two reasons. First, it proves that there exists a function $x^*(\delta)$
which can be used as a solution for every $\delta$.
Moreover, it states that this function can be polynomial. It is not worth detailing
the proof which can
be found in~\cite{bliman04} but it is useful, for the sequel, to recall an outline.\\
~\\
{\bf Outline of the proof of Theorem~\ref{th:Bliman}.} The proof consists in four steps:
\begin{itemize}
 \item[{\bf i)}] It starts with noting that,
if for any $\delta\in K$, there exists $x(\delta)$
such that $G(x(\delta),\delta)>0$ holds, then there exists $\alpha>0$ such that
for any $\delta\in K$, the inequality $G(x,\delta)\geq 2\alpha I$ has a nonempty set of solutions.
It leads to define the set-valued map $F$ which, for any $\delta\in K$, 
associates the nonempty closed convex set
\begin{equation}
F(\delta)=\left\{x\in\R^p:G(x,\delta)\geq \alpha I\right\}.\\
\end{equation}
~\\
(Note that any $x\in F(\delta)$ satisfies ~$G(x,\delta)>0$.)
\item[{\bf ii)}] Then the lower semicontinuity of $F$ is proved.
\item[{\bf iii)}] Michael's selection theorem~\cite{michael} is invoked to prove that, from $F$, there exists a continuous selection $f$
mapping elements of $K$ into elements of $\R^p$ and 
such that, for any $\delta \in K$, $G(f(\delta),\delta)\geq \alpha I$.
\item[{\bf iv)}] The last step applies Weierstrass approximation theorem~\cite{dieudonne}
to prove that each real-valued entry of $f$ is the limit of a sequence of polynomials,
in the sense of uniform convergence.
\end{itemize}
Our purpose leads us to discuss the extension of Theorem~\ref{th:Bliman} where $K$ is a compact subset of $\C^m$
and where each $G_i(\delta)$ belongs to $\C^{n\times n}$.
The three first steps of Bliman's proof can be preserved. All the discussion
concerns the fourth step. We will show that the extension to the complex case is not so obvious and deserves a little attention.

\section{The complex case}
In this section, we discuss to what extend complex parameters can be considered in~Theorem~\ref{th:Bliman}.
But before doing so, we remind the reader of the needed mathematical theorems.

\begin{theorem}
 {\em (Weierstrass approximation theorem~\cite{rudin86})} Suppose that $f$ is a continuous real-valued
function defined on the real range $[a,b]$.  Then
there exists a sequence of real-valued polynomial functions which uniformly converges to $f$.
\end{theorem}
In other words, it is possible to approximate a continuous real-valued function $f$ by a real-valued polynomial
function with any desired accuracy, provided that the degree of the underlying polynomial is large anough.
Bliman actually applies an extension which is the Stone-Weierstrass Theorem.

\begin{theorem}
\label{th:Stone}
{\em (Stone-Weierstrass theorem~\cite{stone48,dieudonne,rudin86})}
Let $K$ be a compact Hausdorff space and let ${\mathcal{P}}$ be a subalgebra of
the set ${\mathcal{C}}(K,\R)$ of continuous functions from $K$ to $\R$,
then ${\mathcal{P}}$ is dense in ${\mathcal{C}}(K,\R)$ if and only if it includes a non
zero constant function and separates the points.
\end{theorem}
Indeed, a subset $K\in\R^m$ is a compact Hausdorff space and the set of polynomial functions with arguments
$\delta$ taken in $K$ is a subalgebra of ${\mathcal{C}}(K,\R)$. It would also be possible to consider $K$ as a subset of
$\C^m$. However, we also want to consider matrices $G(\delta)$ in $\C^{n\times n}$ and therefore
${\mathcal{C}}(K,\C)$, the set of continuous functions from $K$ to $\C$,
should be involved in Theorem~\ref{th:Stone} instead of ${\mathcal{C}}(K,\R)$.
But such a substitution in the writing of the theorem is not possible. The extension is not so simple
but it exists and it is due to Glimm~\cite{glimm60}.
We here propose a simplified version which is especially taylored
for our needs.

\begin{theorem}
\label{th:Glimm}
{\em (Complex version of Stone-Weierstrass theorem~\cite{glimm60})}
Let $K$ be a compact Hausdorff space and let
${\mathcal{C}}(K,\C)$ be the set of continuous functions from $K$ to $\C$.
Also let  ${\mathcal{P}}$ be a subalgebra of ${\mathcal{C}}(K,\C)$ inluding a non
zero constant function and which separates the points.
Then ${\mathcal{P}}$ is dense in ${\mathcal{C}}(K,\C)$ if and only if it is
a C*-algebra with '*' being the relation which associates two complex numbers by conjugation.
\end{theorem}
This theorem allows the consideration of complex functions but also adds an additional property
to be checked by ${\mathcal{P}}$. This is a key issue in our motivation to write this paper
since the set of polynomial functions from $K$ to $\C$ is not stable under conjugation.
Therefore, we have to carefully choose the set ${\mathcal{P}}$ to make this theorem
valid.
\\
Before stating our result, we also recall the next theorem, which is another
avatar of Weierstrass approximation theorem.

\begin{theorem}
\label{th:mergelyan62}
{\em (Mergelyan's theorem~\cite{mergelyan62,rudin86})}
Let $K$ be a compact subset of $\C$ such that $\C\backslash{K}$ is connected.
Any continuous function $f:K\rightarrow \C$ which is holomorphic on its interior $\textup{int}(K)$ can be uniformly
approximated on $K$ by a polynomial function.
\end{theorem}

Bliman's theorem (Theorem~\ref{th:Bliman}) relies on the two first theorems recalled in this section.
With the two other ones, we are now able to state the next theorem which is a key result
before proving Theorem~\ref{th:th2D}.

\begin{theorem}
\label{th:ours}
Let the following expression be defined:

{\scriptsize
\begin{equation}
\label{eq:LMIcomplexe}
 G(x,\delta)=G_0(\delta)+\displaystyle\sum_{i=1}^p x_iG_i(\delta), \delta\in{K},x=\left[\begin{array}{c}
                                 x_1\\
\vdots\\
x_p
\end{array}\right]\in\C^p,
\end{equation}}
where $K$ is a compact subset of $\C^m$
and $G_0$, $G_1$, ..., $G_p$ are continuous mappings defined in $K$ and taking values in
the set of Hermitian matrices of $\C^{n\times n}$.
Also define $K'\subseteq\C^m$ as the image of $K$ under conjugation.
If for any $\delta\in K$, there exists $x(\delta)$
such that $G(x(\delta),\delta)>0$, then there exists $x^*:K\times K'\rightarrow \C^p$,
a polynomial function with respect to $\delta$ and $\delta'$
such that, for any $\delta\in K$, $G(x^*(\delta,\delta'),\delta)>0$.
Moreover, if $m=1$, if $\C\backslash K$ is connected, and if $\textup{int}(K)=\emptyset$,
there exists $x^{\bullet}:K \rightarrow \C^p$, a polynomial function
with respect to $\delta$, such that, for any $\delta$ in $K$, $G(x^{\bullet}(\delta),\delta)>0$.
\end{theorem}

{\bf Proof.} Once again, we remind the reader that the three first steps of the proof
of Theorem~\ref{th:Bliman}~\cite{bliman04} can be reproduced {\em mutatis mutandis}. Only the fourth step
deserves attention. Once we know that there exists $p$ continuous functions from $K$ to $\C$
which are the entries of a solution $x(\delta)$ such that $G(x(\delta),\delta)\geq \alpha I$
for some $\alpha>0$, it remains to prove that, for any given complex entry of $f$,
there exists an approximation.\\
Define the set ${\mathcal{P}}$ of functions of $\delta$ defined through the expressions

\begin{equation}
\label{eq:poly2var}
f(\delta)=\displaystyle\sum_{i=0}^{\mu}\sum_{j=0}^{\mu}c_{ij}\delta^i{\delta'}^j,\quad c_{ij}\in\C.
\end{equation}
These are polynomial functions w.r.t. both $\delta$ and $\delta'$.
It is quite easy to check that ${\mathcal{P}}$ is a subalgebra of ${\mathcal{C}}(K,\C)$
which separates the points and contains non zero constant function.
But in addition, ${\mathcal{P}}$ is closed under conjugation
so it is a C*-algebra.
Therefore, Theorem~\ref{th:Glimm} can be invoked to claim that each entry $x_i(\delta)$,
$i=1,...,p$, of $x(\delta)$ is the limit
of a sequence of polynomial functions complying with~(\ref{eq:poly2var}). Let $x^*(\delta,\delta')$ denote
(for short $x^*(\delta)$)
the $\C^p$-valued function resulting from the concatenation of these limits. Then it comes

\begin{equation}
 \forall \delta\in K, G(x^*(\delta),\delta)\geq \alpha I \Rightarrow \forall
\delta\in K,  G(x^*(\delta),\delta)>0.
\end{equation}
This proves the first part of the theorem.\\
To prove the remainder, it can be noticed that if $\textup{int}(K)=\emptyset$,
the elements of ${\mathcal{P}}$ do not need to be
holomorphic on $K$ but just continuous to invoke Theorem~\ref{th:mergelyan62}.
Therefore, each entry $x^*_i(\delta)$,
$i=1,...,p$ of $x^*(\delta)$ can be uniformly approximated on $K$ by a polynomial
function

\begin{equation}
 x_i^{\bullet}(\delta)=\displaystyle\sum_{i=0}^{\nu}b_i\delta^i,\quad b_i\in\C.
\end{equation}

Let denote $x^{\bullet}(\delta)$
the $\C^p$-valued function resulting from the concatenation of the functions
$x_i^{\bullet}(\delta)$. It follows

\begin{equation}
\forall \delta\in K,\, G(x^{\bullet}(\delta),\delta)>0.\qquad \Box
\end{equation}
The previous theorem deserves several comments. The first statement in Theorem~\ref{th:ours}
is what will really be exploited in the sequel. It shows that the solution to complex parameter-dependent
LMIs can be assumed polynomial with respect to both the parameter vector $\delta$ and its conjugate $\delta'$, with no loss of generality.
This is of course a very interesting consequence of Glimm's theorem but the research for computationally tractable relaxations
of the original LMI conditions often suggests that the solutions should rather be polynomial
or rational with the entries of $\delta$ (not those of $\delta'$, see ~\cite{nDS13,nousTAC14,nDS15}). This additional step on the way to useful approximations
is not so easy to make. Mathematical literature dedicated to approximations usually evokes Mergelyan's theorem (Theorem~\ref{th:mergelyan62})
as the most significant and advanced step but it can be seen that it is restricted to $\delta\in\C$.
Furthermore, $\delta$ should belong to a subset which should satisfy some special properties and these properties
are not encountered in many practical problems, including the present one. For this reason, we claim that the research for efficient approximations
of solutions to complex LMIs with one or several complex parameters, through simple functions, is a widely open problem. More precisely,
consider $G(x,\delta)$ expressed as in~(\ref{eq:LMIreelle}) or~(\ref{eq:LMIcomplexe}) and the underlying function
$g(x,\delta,y)=y'G(x,\delta)y$ which should be non negative: If there are various tools to numerically approximate
solutions $x(\delta,y)$ when $\delta\in\R^m$ and $y\in\R^n$, the existence of such approximations when $\delta\in\C^m$
and $y\in\C^n$ is not as obvious as in the real case. Nevertheless, the very special case addressed in the present contribution
offers a possibility to circumvent the obstacle.\\
If however one wants to enlarge the possibilities of application of Theorem~\ref{th:Glimm} to $m\geq 1$, the following corollary
might be interesting.

\begin{corollary}
\label{cor:corours}
Let an expression be defined by~(\ref{eq:LMIcomplexe})
where $K$ is a compact subset of $\C^m$
and $G_0$, $G_1$, ..., $G_p$ are continuous mappings defined in $K$ and taking values in
the set of Hermitian matrices of $\C^{n\times n}$.
Also define $K'\subseteq\C^m$ as the image of $K$ under conjugation.
If for any $\delta\in K$, there exists $x(\delta)$
such that $G(x(\delta),\delta)>0$, then there exists $x^*:K\times K'\rightarrow \C^p$,
a rational function with respect to $\delta=\left[\begin{array}{ccc}\delta_1' & \dots &\delta_m'\end{array}\right]'$
and $\delta'$, which complies with implicit LFR (ILFR)-based form
\begin{equation}
\label{eq:implicitLFR}
x^*(\delta)=D^*+C^*(E^*-\overline{\Delta}A^*)^{-1}(\overline{\Delta}B^*-F^*),
\end{equation}
where $A^*$, $B^*$, $C^*$, $D^*$, $E^*$ and $F^*$ are either real or complex matrices of appropriate dimensions,
where $\overline{\Delta}=\Delta\oplus\Delta'$ ($\oplus$ denoting ``blocdiag''), where $\Delta=\oplus_{i=1}^m\delta_i I_{q_i}$,
and where $q_i\geq 0$, $i=1,...,m$,
such that, for any $\delta\in K$, $G(x^*(\delta,\delta'),\delta)>0$.
\end{corollary}

{\bf Proof.} (Outline) The proof is more or less the same as the first part of the proof of Theorem~\ref{th:ours}
by noting that when associating two ILFRs of the form~(\ref{eq:implicitLFR}) (by multiplication or addition),
ones gets another ILFR w.r.t. $\overline{\Delta}_1\oplus\overline{\Delta}_2$
where $\overline{\Delta_j}=\Delta_j\oplus\Delta_j'$ and $\Delta_j=\oplus_{i=1}^m\delta_i I_{q_{j_i}}$ (note that the only difference
between $\Delta_1$ and $\Delta_2$ are the dimensions $q_{1_i}$ and $q_{2_i}$). It is well known that the rows and and the columns of
the matrices involved in the ILFR can be re-arranged while preserving the expression and so as it becomes an ILFR
w.r.t. $\overline{\Delta}$ as defined in the corollary. The consequence is that the set of these ILFRs has clearly the properties
of a subalgebra. By applying the conjugation on such an ILFR, another ILFR is obtained, w.r.t. $\overline{\Delta}'=\Delta'\oplus\Delta$.
By again invoking the possible re-arrangement on the rows and the columns of the involved matrices, this ILFR can be written
as an ILFR w.r.t. to $\overline{\Delta}$. Therefore, the set is a C$^*$-algebra and Theorem~\ref{th:Glimm} can be applied.~$\Box$\\
~\\
This corollary will be applied in the next section

\section{Stability of the nD Roesser models}

Now, we come back to our main motivation and prove Theorem~\ref{th:th2D}.
We also propose other perspectives induced by Corollary~\ref{cor:corours}.

\subsection{2D Roesser models}
{\bf Proof of Theorem~\ref{th:th2D}.} Sufficiency holds by virtue of Lemma~\ref{th:intro}.
The interest lies in the second statement and in the assumed structure of $P(\delta)$.
First note that if $P(\delta)$ is solution, so is $P'(\delta)=P(\delta)$ and thus
$P^H(\delta)$. Therefore, there exists $Q(\delta)$ (not necessarily Hermitian) such that,\\
~\\
$\forall\delta\in\partial\cD_2\,(\textup{or} \,\,\forall\delta\in\cD_2^{\diamond}),$
\begin{equation}
\label{eq:steinbis}
\left\{\begin{array}{l}
{\bf M}(\delta)=r_{1_{00}}Q^H(\delta)+(r_{1_{10}}M'(\delta)Q^H(\delta))^H+r_{1_{11}}M'(\delta)Q^H(\delta)M(\delta)<0\\
Q^H(\delta)>0
\end{array}\right.
\end{equation}
Such a system is an LMI system parametered by $\delta$.
Following the notations of the previous sections, it
can be rewritten as $G(x(\delta),\delta)>0$ where $G(\delta)$ is deduced
from $M(\delta)$ and $x(\delta)$ contains the decision variables {\em i.e.}
the entries of $Q(\delta)$~\cite{bgfb94}. The parameter $\delta$
describes the unit circle or the extended imaginary axis which is a compact set $K$.
Matrix $M(\delta)$ is rational w.r.t. $\delta$ and, since $A_{22}$ is stable,
$M(\delta)$ is continuous in $K$ (it has no pole inside $K$). So matrices $G_i(\delta)$ are also continuous.
Therefore, by virtue of Theorem~\ref{th:ours}, $Q(\delta)$ can
be chosen polynomial w.r.t. $\delta$ and $\delta'$
without loss of generality provided that the degree
of the polynomial is large enough. $Q(\delta)$ can then be written

\begin{equation}
\label{eq:Qd_premiere}
 Q(\delta)=\displaystyle\sum_{k=1}^{\eta}\displaystyle\sum_{l=1}^{\gamma}Q_{kl}\delta^k{\delta'}^l,
 \,\,Q_{kl}\in\C^{n\times n}.
\end{equation}
If $\partial\cD_2$ is the unit circle then $\delta \delta'=1$, so each monomial
of the form $\delta^k{\delta'}^l$ equals either $\delta^i$ or ${\delta'}^i$ with
$i=\textup{max}(k,l)-\textup{min}(k,l)$. Therefore $Q(\delta)$ complies with
\begin{equation}
 Q(\delta)=\displaystyle\sum_{i=1}^{\nu}X_{i}\delta^i+Y_{i}{\delta'}^i,
 \,\,(X_{i},Y_{j})\in\{\C^{n\times n}\}^2.
\end{equation}
Since $P(\delta)=Q^H(\delta)$, it can be written as in (\ref{eq:Pdedelta}) with
$P_i=X_{i}+Y_{i}'$.\\
If $\partial\cD_2$ is the imaginary axis, then $\delta'=-\delta$ so~(\ref{eq:Qd_premiere}) obviously
leads to~(\ref{eq:Pdedelta}). Then the proof is completed for $\partial\cD_2$, not for $\cD_2^{\diamond}$.
For that, assume that $P(\delta)$ has the required form~(\ref{eq:Pdedelta}) and is solution for $\delta\in\partial\cD_2$.
Also consider the functions $f_i(\delta)$ defined by

\begin{equation}
f_i(\delta)={\lambda_i}^{-1}(P(\delta)\oplus(-{\bf M}(\delta))),\quad i=1,\dots, 2k_1,
\end{equation}
over $\cD_2^{\diamond}$, where $\lambda_i(.)$ are the eigenvalue of the matrix argument (which are real since the argument
is Hermitian). Since $P(\delta)\oplus(-{\bf M}(\delta))$ is Hermitian,
each eigenvalue is analytic w.r.t. $\delta$.
By invoking the maximum modulus principle~\cite{rudin86}, the functions $f_i(\delta)$
reach their maximum modulus over $\partial\cD_2$. Since from the first part of the proof, $P(\delta)\oplus(-{\bf M}(\delta))$
does not become singular over $\partial\cD_2$, the various $f_i(\delta)$ are bounded over $\partial\cD_2$.
Thus, they are {\em a fortiori} bounded over $\cD_2^{\diamond}$ (they cannot tend to infinity).
Therefore $P(\delta)\oplus(-{\bf M}(\delta))$ never becomes
singular over $\cD_2^{\diamond}$. Consequently, $P(\delta)$ remains positive definite, ${\bf M}(\delta)$
remains negative definite and thus
$P(\delta)>0$ remains a solution
to~(\ref{eq:stein}) over $\cD_2^{\diamond}$.~$\Box$

\begin{remark}
When $\cD_2$ is the open left half complex plane, it is explained in~\cite{nousTAC14} that another admissible expression
of $P(\delta)$ (which could be justified through quite similar arguments) is

\begin{equation}
\label{eq:polylyapbis}
P_G(\delta)=\displaystyle\left(\displaystyle\sum_{h=0}^{\alpha}Q_h\left(\displaystyle\frac{\delta}{1+\delta}\right)^h\right)^H,\quad Q_h\in\C^{k_1\times k_1}.
\end{equation}
Indeed, following the discussion of section~\ref{sec:motiv}, the LFRs involved in the S-procedure used to obtain
non parameter-dependent LMIs would not be well-posed with~(\ref{eq:Pdedelta}) whereas thery are with~(\ref{eq:polylyapbis}).
\end{remark}

It must be noticed that this result generalizes the result obtained, through a different approach, in~\cite{Bliman2002}
for the discrete case. Here, all the cases are embedded in the same framework.

\subsection{Insights into $n$D Roesser models}

Consider a generalization of model~(\ref{eq:Roesser}) to dimension $n$:

\begin{equation}
\label{eq:Roesser_n}
\left[\begin{array}{c}
       {\bf q}_1x^1(j_1,...,j_n)\\
       \vdots\\
{\bf q}_nx_2(j_1,n...,j_n)
      \end{array}\right]=\left[\begin{array}{ccc}
A_{11} & \dots & A_{12}\\
\vdots & \ddots & \vdots\\
A_{n1} & \dots & A_{nn}
\end{array}\right]\left[\begin{array}{c}
       x_1(j_1,...,j_n)\\
       \vdots\\
x_n(j_1,\dots,j_n)
      \end{array}\right],
\end{equation}
where ${\bf q}_i$ denotes either shift operator or derivative along dimension $i$. In many applications
of $n$D models, one dimension (most probably the time) is not considered on an equal footing and one is mainly
interest in the stability w.r.t. that dimension. Let us consider with no loss of generality
that this dimension is the first one. Then consider the transform $\cX: x_i(j_1,j_2,,...,j_n)\mapsto X_i(j_1,\lambda_2,...,\lambda_n)$
where $\cX$ actually corresponds to the Laplace transform for the continuous dimensions and the ${\cal Z}$ transform
on the discrete dimension. If we define a ``mixed'' ${\cal L}_2-{\mathit l}_2$-norm of $X_i$, denoted by
$||X_i||_2$, then the dynamics of $x_1$ is related to the behaviour of $||X_1||_2$ under the constraint

\begin{equation}
{\bf q}_1X_1(j_1,\lambda_2,...,\lambda_n)=\underbrace{\mathop{D_{M}+C_M
\left(I-
 \tilde{\Delta}A_M\right)^{-1}\tilde{\Delta}B_M}}_{M(\delta)} X_1(j_1,\lambda_2,...,\lambda_n),
\end{equation}
where $j_1$ is a continuous or discrete time, $\lambda_i$, $i=2,...,n$ describes either the unit circle or the imaginary axis,
where $\delta=\left[\begin{array}{ccc}(\delta_1=\frac{1}{\lambda_2})' & \dots & (\delta_{n-1}=\frac{1}{\lambda_n})'\end{array}\right]'$,
where $\tilde{\Delta}=\oplus_{i=1}^n\delta_{n-1}$, and where

\begin{equation}
\left[\begin{array}{c|c}
       A_M & B_M\\
       \hline
       C_M & D_M
      \end{array}\right]=\left[\begin{array}{ccc|c}
       A_{22} & \dots & A_{2n} & A_{21}\\
\vdots & \ddots & \vdots & \vdots\\
A_{n1} & \dots & A_{nn} & A_{n1}\\
\hline
A_{12} & \dots & A_{1n} & A_{11}
      \end{array}\right].
\end{equation}
Thus, in the sense of $||.||_2$, the system is asymptotically stable along dimension~$1$ if and only if
$M(\delta)$ has no eigenvalue outside of $\cD_1$ ({\em i.e.} outside of $\overline{\C}^-$ or ${\bf D}$)
for any $\delta_i$ describing either the unit circle or the imaginary axis, depending on the nature of the propagation
along dimension $i+1$. This straightforwardly results in the next theorem.

\begin{theorem}
\label{th:nD}
System~(\ref{eq:Roesser_n}) is {\em time-relevant} stable in the sense of $||.||_2$ if and only if
Inequality~(\ref{eq:stein}) holds for any $\delta\in\partial\cD_2\times ...\times \partial\cD_n$,
where the various $\partial\cD_i$ are defined by~(\ref{eq:defDi}) and by generalizing the constraint~(\ref{eq:choices_q})
to  dimension $n$.
\end{theorem}
Therefore we state the next theorem

\begin{theorem}
\label{th:nDbis}
If system~(\ref{eq:Roesser_n}) is time-relevant stable in the sense of $||.||_2$, then there exists a solution $P^{\bullet}(\delta)$
to Inequality~(\ref{eq:stein}) over $\partial\cD_2\times ...\times \partial\cD_n$
which complies with
the ILFR-based expression
\begin{equation}
\label{eq:Pbullet}
P^{\bullet}(\delta)=D^{\bullet}+C^{\bullet}(E^{\bullet}-\hat{\Delta}A^{\bullet})^{-1}(\hat{\Delta}B^{\bullet}-F^{\bullet}),
\end{equation}
where $\hat{\Delta}=\oplus_{i=1}^m \delta_iI_{g_i}$ and $g_i\geq 0$,  $i=1,...,(m=n-1)$.
\end{theorem}

{\bf Proof.} From Corollary~\ref{cor:corours}, it can be easily assessed that there exists a solution
$P^*(\delta)=D^{*}+C^{*}(E^{*}-\overline{\Delta}A^{*})^{-1}(\overline{\Delta}B^{*}-F^{*})$ with $\overline{\Delta}=\Delta\otimes \Delta'$
$\Delta=\oplus_{i=1}^m\delta_i I_{q_i}$, and $q_i\geq 0$. Since $\delta_i$ belongs to either the imaginary axis
or the unit circle, it can be written $-\delta_i$ or $\delta_i^{-1}$. Those two expressions are also special
instances of ILFRs. Since ILFRs of ILFRs are also ILFRS, then it is possible to rewrite the solution $P^*$
as $P^{\bullet}$ given in~(\ref{eq:Pbullet}).~$\Box$\\
~\\

The question is now that of (possibly exactly) relaxing conditions such as~(\ref{eq:stein}). It cannot be addressed
here unless by unduly increasing the size of the article. However, when $n=2$ ({\em i.e.} $m=1$), Theorem~\ref{th:th2D}
was exactly relaxed in~\cite{nDS13} in the purely discrete case through S-prcocedure. Using a close by different approach,
the mixed case was efficiently tackled in~\cite{ChesiMiddletonTAC14} where SOS-relaxations were prefered to S-procedure.
Moreover, in~\cite{ChesiMiddletonTAC14}, upperbounds on the degree of the polynomial are provided, which is very interesting
from a theoretical point of view. But the bounds are still far larger than the degrees actually required to assess stability,
and therefore might not be easily exploited on practical examples.\\
An exact LMI-based relaxation is provided in~\cite{Ebihara2006-1}. It is very interesting because unlike the other approaches,
it does not require any degree to be guessed. Therefore, instability can be proved by such a condition. Nevertheless,
it is restricted to the discrete case, and considers Fornasini-Marchesini (FM) models instead of Roesser (R) models.
Of course, FM and R models can be converted to each other in the control-free case but it is no longer true in the presence
of a control input. So any attempt to extend the result to design may clearly separate the two models. Moreover, examples
exposed in~\cite{nousTAC14} tend to show that the condition of~\cite{Ebihara2006-1}, unless more direct, may not necessarily be more
attractive than~\cite{ChesiMiddletonTAC14,nousTAC14} in terms of computation time. Furthermore, we would like to insist
that~\cite{nousTAC14} exploits Theorem~\ref{th:th2D} (or more precisely a slightly modified version)
to encompass all the cases (discrete, continuous, and mixed).\\
At last, in~\cite{nDS15}, the first exploitation of Theorem~\ref{th:th2D} is made for control and a nesssary and sufficient
condition for state feedback structural stabilization of 2D discrete models is proposed. This is a major advantage
compared to other approaches. The generalization to continuous and mixed cases is under investigation.\\
For $n>2$ (or $m>1$), the S-procedure-based reasoning used in~\cite{nousTAC14} might be extended but also probably
lead to a conservative LMI condition, hence our focus on the 2D case.\\
Once again, the reader must be aware that the previously evoked results~\cite{nDS13,nousTAC14,nDS15} rely on Theorem~\ref{th:th2D}
whereas it had not been clearly proved yet. The main contribution of the present paper is to fill this gap and to embed this result
into a more general formalism.

\section{Conclusion}
In this note, we proposed a necessary and sufficient condition for the stability of 2D (discrete, continuous or mixed) Roesser models
which consists of an LMI system where the solutions are polynomial w.r.t. a complex parameter describing either the unit circle
or the imaginary axis (or either the closed unit disc or the closed right half-plane). Insights to generalize the results
to higher dimensions were also proposed.
It has to be noticed that
this result is in accordance with what was obtained in the discrete case through a completely different approach in~\cite{Bliman2002}
and that a recent paper proposed a quite similar result in the mixed continuous-discrete-case~\cite{ChesiMiddletonTAC14}.
But the present approach can offer a general method to derive exact relaxations of 2D-stability conditions
for all the cases~\cite{nousTAC14}. At last, but not least, a first exploitation of this result for the synthesis
of control laws was made in~\cite{nDS15}. Therefore, Theorem~\ref{th:th2D} is a keystone of a promising approach.

\bibliographystyle{plain}        

\end{document}